\documentstyle[12pt]{article}

\font\tenmsb=msbm10 scaled \magstep1
\font\sevenmsb=msbm7 scaled \magstep1
\font\fivemsb=msbm5 scaled \magstep1

\newfam\msbfam
\textfont\msbfam=\tenmsb
\scriptfont\msbfam=\sevenmsb
\scriptscriptfont\msbfam=\fivemsb
\def\Bbb#1{{\fam\msbfam\relax#1}}

\def\R{{\Bbb R}}
\def\C{{\Bbb C}}
\def\N{{\Bbb N}}

\def\Subset{\subset\subset}
\def\e{\varepsilon}
\title{
    \hfill\raisebox{0.0ex}[1.ex][3ex]{}\newline
\bf Solutions of algebraic equations \\
with analytic almost periodic coefficients\footnote {This research was 
supported by INTAS-99-00089}}
\author{\bf Britik V.V., Favorov S.Ju.}
\date{ }
\begin{document}
\maketitle

It is well known that solutions $w (z )$
of the equation
$$
      a_{m}(z ) w^{m} + a_{m-1}(z )
      w^{m-1} + \dots + a_{1}(z ) w + a_{0}
      (z ) = 0                               \eqno(1)
$$
often inherit properties of the coefficients $a_{j}(z ),
\phantom{2} j=0,\dots,m,$. As an example,
suppose that these coefficients are almost periodic functions on the
axis, $a_{m}(z )=1$, and the discriminant $D(z
)$ of the polynomial in (1) satisfies the condition
$$
      |D(z )| \geq \gamma >0;     \eqno(2)
$$
then each solution  of (1) is an almost periodic function, too (\cite{c},
\cite{b}).
Nevertheless, one cannot replace condition (2) by the weaker condition
$$
      D(z) \neq 0                           \eqno(3)
$$
even for the equation
$$
      w^2 - a_0(z ) = 0                      \eqno(4)
$$
(\cite {w}). However for analytic almost periodic coefficients
$a_{j}(z ),\phantom{2} j=0,\dots,m$,\, on a strip $S$,
the conditions $a_{m}(z )=1$ and (3) imply that every
continuous solution of (1) is an analytic almost periodic
function on this strip (\cite{f}).

      Note also that one can formulate classical Bohr's
theorem on division of analytic almost periodic functions
(see for example \cite {l})
in the following way: an analytic solution of (1) for
$m=1$ and analytic almost periodic functions $a_{1}(z)$,
$a_{0}(z)$ on a strip is an almost periodic
function on this strip.

      It is natural to consider analytic solutions of (1) with
analytic almost periodic coefficients without any restriction
on the discriminant $D(z )$. We know only one result of
this kind: namely, an analytic solution of (4) with an analytic
almost periodic function $a_{0}(z )$ on a strip  is
almost periodic as well. However, by our opinion, the proof of this result
in \cite {j} is not perfect.

    Recall that a function $f(z)$ is said to be {\it almost
 periodic on the real axis} $\R$ if $f(z )$ belongs to the
closure  of the set of finite exponential sums
$$
\sum{a_{n}e^{i\lambda_{n}z}}, \phantom{e} a_n\in\C, \phantom{o}
\lambda_n\in\R,                                            \eqno(5)
$$
with respect to the topology of uniform convergence on $\R$.
Further, let $S$ be a strip
$\left\{z\in\C: a<{\rm Im}z<b\right\}$
($a$ can be $-\infty$ and $b$ can be $+\infty$). We write $S'\Subset S$
if $S'=\left\{z\in{C}:a'<{\rm Im}z<b'\right\}$,\quad $a<a'<b'<b$.
A function $f(z)$ is said to be {\it analytic almost
periodic on a strip} $S$ if $f(z)$ belongs to the closure
of the set of sums (5) with respect to the topology of uniform
convergence on every substrip  $S'\Subset S$. The equivalent
definitions are the following: the family $\left\{f({z+h}
 )\right\}_{h\in\R}$ is a relative compact set
 with respect to the  topology of uniform convergence on $\R$ (for
almost periodic functions on the axis) or with respect to the
topology of uniform convergence on every substrip $S'\Subset S$
(for analytic almost periodic functions).

By $AP(S)$ we denote the space of all analytic almost periodic
functions on $S$ equipped with the topology of uniform
convergence on every substrip $S'\Subset S$;
the zero set of a function $f\in{AP(S)}$ is denoted by $Z(f)$.

\medskip

      {\bf Theorem 1}. {\it Let $w(z)$ be a
continuous solution of (1) in $S$ and $a_{k}{(z)}
\in{AP(S)}$, $0\leq{k}\leq{m}$. Then $w(z
)\in{AP(S)}$.}

    Proof of this theorem makes use of the following simple lemmas
on roots of polynomials
$$
      Q(w ) =w^m + b_{m-1}w^{m-1} + \dots +b_1w + b_0.  \eqno(6)
$$

\medskip

      {\bf Lemma 1}. {\it For any $N<\infty$, $\varepsilon>0$,
there exists a constant $\nu>0$ depending on $N$ and
$ \varepsilon$ only such that the roots $w_j,\ j=1,\dots,m$, and
$\tilde w_j,\ j=1,\dots,m,$ of any polynomials $Q$, $\tilde{Q}$
of the form (6) with $\max_j|b_j|\leq N,\  \max_j|\tilde b_j|\leq N,\
\max_j |b_j-\tilde{b_j}|\leq{\nu}$,
satisfy, under a suitable numeration, the conditions
$|w_j-\tilde{w_j}|<\varepsilon$, $j=1\dots m$.}

{\bf Proof}. Assume the contrary. Then for some $N<\infty,\ \e_0>0$\
there exist two sequences of polynomials
$$
Q_n(w ) =w^m + b^{(n)}_{m-1}w^{m-1} + \dots + b^{(n)}_0,\quad
\tilde Q_n(w) =w^m +\tilde b^{(n)}_{m-1}w^{m-1} + \dots +\tilde b^{(n)}_0
$$
such that $\max_j|b_j|\leq{N},\ \max_j|\tilde b_j|\leq N,\
\max_j |b_j-\tilde{b_j}|\to 0$ as $n\to\infty$, and
$$
\max_j |w_j^{(n)}-\tilde w_j^{(n)}|>\e_0\quad\eqno(7)
$$
under every numeration of roots $w_j^{(n)},\ \tilde w_j^{(n)},
\ j=1,\dots,m,$ of the polynomials $Q_n$,
$\tilde{Q_n}$ respectively.
Without loss of generality it can be assumed that
$$
b^{(n)}_j\to\overline b_j,\quad \tilde b^{(n)}_j\to\overline b_j,
\quad j=0,\dots,m-1,\quad {\rm as}\ n\to\infty.
$$
Hence the sequences of the polynomials $Q_n,\  \tilde{Q_n}$
converge to the same polynomial
$$
    \overline Q(w)=w^m +\overline b_{m-1}w^{m-1} + \dots +\overline b_0
$$
with respect to the uniform convergence on every compact subset of $\C$.

Let $C_j,\ j=1,\dots,p,\ p\le m,$ be disjoint disks of radius $r<\e_0/2$
with the centers at the roots of the polynomial  $\overline Q(w)$.
By Hurwitz' Theorem, for $n$ large enough all roots of the  polynomials
$Q_n(w)$, $\tilde Q_n(w)$, lie in these disks and a number of roots of
the polynomial $Q_n(w)$ \ in  a disk $C_j$ coincides with
a number of roots of the polynomial $\tilde Q_n(w)$\ in the same disk
for each $j=1,\dots,p$. Therefore there exists a numeration of
roots of the polynomials $Q_n$, $\tilde{Q_n}$ such that (7) is false.
This contradiction proves the lemma.

\medskip

{\bf Lemma 2}. {\it The distance between any two roots
of a polynomial $Q$ of the form (6) with the discriminant
$d(Q)\neq 0$ is greater
than a constant $\tau>0$ depending on $|d(Q)|$ and  $\max_j|b_j|$ only.}

{\bf Proof}. Assume the contrary. Then there exists a sequence
of polynomials
$$
      Q_n(w ) =w^m + b^{(n)}_{m-1}w^{m-1} + \dots + b^{(n)}_0
$$
such that $\max_j|b_j|\le N<\infty,\ |d(Q_n)|\ge\delta>0$, and the distance
between some two roots of a polynomial $Q_n$ tends to zero as $n\to\infty$.
Without loss of generality it can be assumed that
$b^{(n)}_j\to\overline b_j,\ j=0,\dots,m-1,$ as $n\to\infty$; hence
the discriminants $d(Q_n)$ converge to the discriminant $d(\overline Q)$
of the polynomial
$$
    \overline Q(w)=w^m +\overline b_{m-1}w^{m-1} + \dots +\overline b_0
$$
and $d(\overline Q)\neq 0$. Using Lemma 1 for $\overline Q$ and $Q_n$ with
$n$ large enough, we obtain that the distance between some two roots of
the polynomial $\overline Q$ is arbitrary small, i.e., this polynomial
has a multiple root. This contradicts the assertion
$d(\overline Q)\neq 0$. The lemma is proved.

\medskip

      {\bf Proof of Theorem 1}. We may assume that $a_m(z)
\not\equiv{0}$. First let us suppose that the discriminant
$D(z)\not\equiv{0}$. The solution $w(z)$ is bounded on a
neighborhood of any point $z'\in S$; moreover $w(z)$ is analytic
at any point $z'\in S$ such that $a_m(z')\not=0$, $D(z')\not={0}$.
Since zeros of $a_m(z)$ and $D(z)$ are isolated,  $w(z)$
is analytic on $S$.

Let us show that for an arbitrary sequence
$\{h_n\}\subset\R$ there exists a subsequence $\{h_{n'}\}$ such
that the functions
$w(z+h_{n'})$ form a Cauchy sequence in the space $AP(S)$.
   It is sufficient to check that these functions converge
uniformly on each substrip $S_0\Subset S_1\Subset S$. We may assume that
the  functions  $a_j(z+h_{n})$ converge to functions
$\overline{a_j}(z)$ in the space $AP(S)$ for $n\to{\infty}$  and  each
$j=1,\dots,m$; then the functions $D(z+h_{n})$ converge in  this
space to the discriminant $\overline D(z) $ of the right  hand
side  of the  equation
 $$
  \overline a_{m}(z)w^{m}+\dots+\overline a_{1}(z)w+\overline
       a_{0}(z)w=0                                        \eqno(8)
$$
  Let $U_r$ be the $r$-neighborhood of the set
 $[Z(\overline a_{m})\bigcup{Z(\overline{D})}]\bigcap{S_1}$.  We claim  that
for sufficiently small $r$ there exist closed rectangles $\Pi_l$
such that
 $S_0\subset{\bigcup \limits_{l} \Pi_{l}}\subset {S_1}$
and $\partial \Pi_{l}$ disjoint with $U_r$ for all $l\in\N$.

    Since the functions $\overline D(z)$ and $\overline a_{m}(z)$ belong to
$AP(S)$,  the numbers of their zeros
inside a rectangle
$\{z\in S_1:|{\rm Re}z-t|<1\}$
are bounded by a number $K$ independent of $t\in\R$
(see \cite{l}). Hence for $r<1/4K$ and for all $t\in\R$ there
exists $c_{t}\in\R$ such that $ |c_{t}-t|<1 $ and the straight line
${\rm Re}z=t$ does not intersect the set $U_r$. Then there exists
a sequence of rectangles $\{ z\in{S_1} : c_l\le{\rm Re}z\le{c_{l}'}\}$
overlapping the strip $S_1$ whose lateral sides
are disjoint with $U_r$. Furthermore, suppose
$r<(8K)^{-1} \inf\{|z-z'|: z\in{S_0},\phantom{o} z'\notin{S_1}\}$,
then there exist segments
$\{z:{\rm Im}z=d_l,c_l\le{{\rm Re} z}\le{c_{l}'}\}
\subset {S_1\setminus{S_0}}$
disjoint with $U_r$ as well. Thus the rectangles
$\Pi_l=\{z: c_l\le{{\rm Re}z}\le{c_l'}, \,d_l\le{{\rm Im}z}\le{d_l'}\}$
with  suitable $d_l$, $d_l'$ are just required.

   It follows from properties of analytic almost periodic functions
(see \cite{l}) that $|\overline{D}(z)|\ge \eta$,\
$|\overline{a_{m}}(z)|\ge \eta$
for $z\in {S_1\setminus{U_r}}$, where $\eta$ is a strictly positive
constant. Hence for $n\ge N$ and $z\in {S_1\setminus{U_r}}$, we have
$|D(z+h_n)|\ge \eta,\,|a_m(z+h_n)|\ge \eta$.
Besides, the functions $a_j(z), 0\le j\le m-1,$ are uniformly
bounded on $S_1$. Applying Lemma 2, we get that the distance
between any two roots of the polynomial
$$
     Q_{n} (w ) = w^m + \frac{a_{m-1}(z +
      h_{n} )}{a_{m}(z + h_{n} )} w^{m-1}
      + \dots + \frac{a_{0}(z + h_{n} )}{a_{m}
      (z + h_{n} )}
$$
is greater than $\tau >0$. Note that the constant $\tau $ is the
same for all $z\in {S_1\setminus{U_r}}$, $n\ge N$. Further, the
functions
$$
{a_{j}(z + h_n)\over a_m(z+h_n)}
$$
form a Cauchy sequence with respect to the uniform convergence on the set
$S_1\setminus{U_r}$ for every $j=0,\dots,m-1$. This yields that the
polynomials $Q_{n} (w )$, $Q_{k} (w )$ satisfy
the conditions of Lemma 1 with $\varepsilon =\tau/3$ and
$n,k \ge N_1(\varepsilon)$ for all $z\in {S_1\setminus{U_r}}$.
Hence for every fixed $z\in {S_1\setminus{U_r}}$ there exists a
solution $\tilde w(z)$ of the equation $Q_n(w )=0$
such that
$$
|w(z+h_k)-\tilde w(z)|\le  \varepsilon.
$$

   Now we have two possibilities for each $z\in {S_1\setminus{U_r}}$ :
either
$\tilde w(z)=w(z+h_n)$ and
$$
      |w(z+h_k)-w(z+h_n)|\le \frac{\tau}{3}, \eqno(9)
$$
or  $|\tilde w(z)- w(z+h_n)|\geq\tau$  and
$$
      |w(z+h_n)-w(z+h_k)|\geq |w(z+h_n)-\tilde w(z)|-
      |\tilde w(z)-w(z+h_k)|\ge \frac{2}{3} \tau.
$$

    Fix an arbitrary point $z_0\in  \partial \Pi_{1}$. The coefficients
of the polynomials $ Q_{n} (w )$ are bounded at this point,
therefore the sequence $w(z_0+h_n)$ is also bounded. Without loss of
generality it can be assumed that this sequence converges, hence
inequality (9) is true for $z=z_0$. Since the set
$\bigcup \limits_{l} \partial \Pi_{l}$ is connected, we see
that (9) holds on this set. Using the Maximum Principle, we obtain that
(9) is true for all $z\in\bigcup_l \Pi_l \supset S_0$.
Hence we have $\tilde w(z)= w(z+h_n)$ for all $z\in S_0$. Thus
the functions $w(z+h_n)$ form a Cauchy sequence with respect to the
uniform  convergence on $S_0$ and $w(z)$ is an almost periodic function
on $S$.

    If the discriminant of the polynomial
$P(w)=a_m(z)w^{m}+\dots+a_1(z)w+a_0(z)w$
is zero, then the equations $P(w)=0$ and $P'(w)=ma_m(z)w^{m}+\dots
+a_1(z)=0 $ have a common solution for each fixed $z\in S$. Using the
Euclid algorithm, we get
$$
      P(w)=Q(w)R(w),\phantom {dtydtdtyd} P'(w)=T(w)R(w),
$$
      where the coefficients of $Q(w)$, $T(w)$, $R(w)$ lie in the
quotient field of $AP(S)$. Besides, if $w(z)$ is a solution of (1)
for fixed  $z\in S$, then $w(z)$  is an ordinary solution of the
equation $Q(w)=0$ whenever all the coefficients of $Q(w)$ are finite
at this point $z$. Multiplying $Q(w)$ by a suitable function
from $AP(S)$, we obtain a polynomial $\tilde Q(w)$ with the
coefficients from $AP(S)$  such that $w(z)$  is an ordinary
solution of the equation $\tilde Q(w)=0$ for all $z\in S$ outside
of some discrete set. Hence  the discriminant of $\tilde Q(w)$ does
not vanish and we can use the previous result. The theorem is proved.

\medskip

      {\bf Theorem 2}. {\it Suppose $w(z)$ is a meromorphic
solution of (1) with $a_j(z)\in AP(S),\ j=0,\dots,m,$ and
$$
    {\rm card} \{ z \in S^{\prime} : | {\rm Re} \phantom{.}
      z | < t, \phantom{1} w(z) = \infty \} =
      o(t) \phantom{12} as \phantom{12} t \to
      \infty
$$
for each  $S'\Subset S_1$. Then $w(z)\in AP(S)$.
}

   {\bf Proof.} Let $S_0\Subset S_1\Subset S$.
It can be easily seen that for all $t\in\R$ there exists a rectangle
$\{ z\in S_1 : |{\rm Re}z-h|<t \}$ without poles of $w(z)$. Hence there
exists a sequence of rectangles $\{ z\in S_1 : |{\rm Re}z-h_n|<t_n \},
\phantom{12} t_n \to \infty $, without poles of $w(z)$.
We may assume $a_m(z)\not\equiv 0$, $D(z)\not\equiv 0$ and the sequences
of the functions $a_j(z+h_n), \, j=0,\dots,m,$ \ $D(z+h_n)$
converge in the space $AP(S)$ to functions $\overline a_j(z)$,
$\overline D(z)$ respectively. Note
that all poles of $w(z)$ lie in the set $Z(a_m)$. Applying the
arguments of Theorem 1, we obtain that the sequence  $w(z+h_n)$
converges uniformly on the set $\bigcup\limits_l\partial\Pi_l$.

     Let $\overline w(z)$ be the limit of the sequence. Since every
rectangle $\Pi_{l}$ lies inside the set $\{ z\in S_1 : |{\rm }z|<t_n \}$
for $n \ge n(l),$ we see that the functions  $w(z+h_n)$ converge on
$\Pi_{l}$ to an analytic function, therefore $\overline w(z)$ is
analytic on $S_0$. Now Theorem 1 implies that $\overline w(z)$ is an
almost periodic solution of (8).

     Furthermore,
$$
      \sup_{S^{\prime}} |\overline{a_{j}} ( z -
      h_{n} ) - a_{j} (z)| =
      \sup_{S^{\prime}} |a_{j} ( z +
      h_{n} ) - \overline{a_{j}} (z)
      | \rightarrow 0 \phantom{12} {\rm as} \phantom{12}
      n \rightarrow \infty
$$
for each  $S'\Subset S$ and $j=0,\dots,m.$
Applying the above arguments, we see that the sequence of the
functions $\overline w(z-h_n)$ converges in the space $AP(S_0)$
to an analytic solution $\overline{\overline{w}}(z)$ of (1).

   Since $S_0$ is an arbitrary substrip of $S$, we only need
to prove that $\overline{\overline w}(z)=w(z)$. Assume the contrary.
Let $S'\Subset S_0$ be an arbitrary substrip,
$\tilde U_r$ be the $r-$neighborhood of the set
 $[Z(a_{m})\cup Z(D)]\cap S'$. Applying
the above arguments and Lemma 2 we see that
$$
|\overline{\overline{w}}(z)-w(z)|\ge\tau>0 \phantom{..x}{\rm for\
all}\phantom{.}z\in S'\setminus \tilde U_r   \eqno(10)
$$
with certain $\tau>0$. On the other hand, we have
$$
|w(z+h_n)-\overline{w}(z)|\le\tau/3 \phantom{x..}{\rm for}
\phantom{.}n\ge n(\tau),\phantom{.} z\in\bigcup\limits_l\partial\Pi_l.
$$
Therefore, using the uniform convergence of $\overline{w}(z-h_n)$
to $\overline{\overline{w}}(z)$ on $S'$, we obtain

$$
|\overline{\overline{w}}(z)-w(z)|\leq
      |\overline{\overline{w}}(z)-\overline{w}(z-h_n)|+
      |\overline{w}(z-h_n)-w(z)|\leq \frac{2}{3}\tau
   \eqno(11)
$$
for $z\in\bigcup\limits_l\partial\Pi_l-h_n$ and  sufficiently
large $n$. Arguing as in the proof of Theorem 1, we can see that 
if $r$ is small enough, then every
vertical segment $\{ z\in S': \ {\rm Re} z=t\} $ has common points
with $S'\setminus\tilde U_r$. Thus inequalities
(10) and (11) are simultaneously fulfilled on the nonempty set. 
This contradiction completes the proof.

     The authors are grateful to the late Professor L.I. Ronkin
who had called their attention to the problem considered in this
paper.

\bigskip
{\it Harkiv National university

md.Svobody 4, Harkiv, 61077, Ukraine.

\medskip
e-mail: favorov@ilt.kharkov.ua}

\vskip 1.5cm

Britik V.V.,Favorov S.Ju. {\bf Solutions of algebraic equations with
analytic almost periodic coefficients.}
We prove that continuous or meromorphic, with a small number of poles,
solutions of algebraic equations with the analytic almost periodic
coefficients are almost periodic, too.

\end{document}